\documentclass[12pt]{amsart}
\usepackage{amscd,amssymb,verbatim}
\usepackage{amssymb,amsmath,amsthm,epsfig}

\def\cnum#1{\bigcirc\kern -8pt#1}
\newcommand{\n}{\noindent}

\theoremstyle{plain}
\newtheorem{thm}{Theorem}[section]

\newtheorem{cor}[thm]{Corollary}
\newtheorem{prop}[thm]{Proposition}
\newtheorem{lem}[thm]{Lemma}
\newtheorem{defn}[thm]{Definition}
\newtheorem{rem}[thm]{Remark}

\newcommand{\N}{\mathbb{N}}

\newcommand{\mnorm}[1]{\left\vert\kern-0.9pt\left\vert\kern-0.9pt\left\vert #1
    \right\vert\kern-0.9pt\right\vert\kern-0.9pt\right\vert}

\begin{document}

\title{Embedding $\ell_{\infty}$ into the space of all Operators on Certain Banach Spaces}

\author{G. Androulakis, K. Beanland$^*$, S.J. Dilworth, F. Sanacory$^*$}
\date{}
\thanks{$^*$The present paper is part of the Ph.D theses of the second
and fourth author which are prepared at the University of South Carolina under the
supervision of the first author}

\begin{abstract}
We give sufficient conditions on a Banach space $X$ which ensure
that $\ell_{\infty}$ embeds in $\mathcal{L}(X)$, the space of all
operators on $X$. We say
that a basic sequence $(e_n)$ is quasisubsymmetric if for any two
increasing sequences $(k_n)$ and $(\ell_n)$  of positive integers
with $k_n \leq \ell_n$ for all $n$, we have that $(e_{k_n})$
dominates $(e_{\ell_n})$. We prove that if a Banach space $X$ has a
seminormalized quasisubsymmetric basis then $\ell_{\infty}$ embeds
in $\mathcal{L}(X)$.
\end{abstract}

\maketitle

\markboth{G. Androulakis, K. Beanland, S.J. Dilworth, F. Sanacory}{Embedding $\ell_{\infty}$
into the space of all Operators on Certain Banach Spaces}

\section{Introduction} \label{sec1}

The famous open problem of whether there exists an infinite dimensional Banach space on
which
every (linear bounded) operator is a compact perturbation of a multiple of the identity,
is attributed to S. Banach (see related papers \cite{G}, \cite{GM}, and \cite{S}).
One of the reasons that this problem has attracted a lot of attention is that if
such a space $X$ exists then by the results of \cite{AS} or \cite{L}, $X$ provides a positive
solution to the Invariant Subspace Problem for Banach spaces, namely every operator on
$X$ has a non-trivial (i.e. different than zero and the whole space) invariant subspace.
Notice that if a space $X$ satisfies the assumptions of the above problem of Banach and
if in addition $X$ has the Approximation Property and separable dual then $\mathcal L (X)$,
the space of all operators on $X$ endowed with the usual operator norm, is separable.
Thus if a reflexive Banach space $X$ with a basis satisfies the assumptions of the above
problem of Banach then $\mathcal L (X)$ is separable. In this paper we provide sufficient
conditions on a Banach space $X$ which imply that $\ell_{\infty}$ embeds in
$\mathcal{L}(X)$ (which we denote by ``$\ell_{\infty} \hookrightarrow \mathcal{L}(X)$'')
hence $\mathcal{L}(X)$ is non-separable. The question whether ${\mathcal L}(X)$ is
non-separable for every infinite dimensional Banach space $X$, is a well known open
problem.

In Section 2 we introduce the notion of a quasisubsymmetric basic
sequence, give examples and prove our main results. We say that a
basic sequence $(e_n)$ is {\em quasisubsymmetric} if for any two
increasing sequences $(k_n)$ and $(\ell_n)$  of positive integers
with $k_n \leq \ell_n$ for all $n$, we have that $(e_{k_n})$
dominates $(e_{\ell_n})$. One of our main results
(Theorem~\ref{th:QSSEmbedding}) is that if a Banach space $X$ has a
seminormalized quasisubsymmetric basis then $\ell_{\infty}$ embeds
in $\mathcal{L}(X)$. Thus a reflexive Banach space $X$ with a
seminormalized quasisubsymmetric basis can not satisfy the
assumptions of the above problem of Banach. Moreover, in
Proposition~\ref{prop:nseparable} we prove that if a Banach space
$X$ has a seminormalized basis which dominates all of its
subsequences then $\mathcal{L}(X)$ is non-separable (we do not know
whether $\ell_{\infty} \hookrightarrow \mathcal{L}(X)$ in this
case). Finally we observe how our results provide sufficient conditions on a 
Banach space $X$ which ensure that $\ell_1$ embeds complementably 
into the space of nuclear operators of $X$.

In the rest of the section we review some known results about embedding $\ell_{\infty}$
into the space of operators. If $X$, $Y$ are Banach spaces then denote by
${\mathcal L}(X,Y)$ the space of all operators from $X$ to $Y$ and ${\mathcal K}(X,Y)$
the space of all compact operators from $X$ to $Y$.

N.J.~Kalton \cite[proof of Theorem 6 (iii) $\Rightarrow$ (ii)]{K}
proved that if $X$ is an infinite dimensional separable Banach space, $Y$ is any Banach
space then:
\begin{itemize}
\item[(a)] if $c_0$ embeds in ${\mathcal L}(X,Y)$ then $\ell_{\infty}$ embeds in
${\mathcal L}(X,Y)$,
\item[(b)] if $c_0$ embeds in $Y$ then $\ell_{\infty}$ embeds in ${\mathcal L}(X,Y)$.
\end{itemize}
Later M. Feder \cite[page 201]{F} noticed that by the result of
B.~Josefson \cite{J} and A.~Nissenzweig \cite{N}, the proof of
Kalton for the above results works without the assumption that $X$
is separable. Notice that the result (a) of Kalton as extended by
Feder generalizes the result of C. Bessaga and A. Pelczynski
\cite{BP} stating that if $c_0$ embeds in $X^*$ then $\ell_{\infty}$
embeds in $X^*$.

A.E.~Tong and D.R.~Wilken \cite{TW} proved that if $X$, $Y$ are Banach spaces, $Y$ has an unconditional basis and there is a noncompact operator in ${\mathcal L}(X,Y)$
then $\ell_{\infty}$ embeds in ${\mathcal L}(X,Y)$. The main focus of the paper of
Tong and Wilken is the well known open problem of whether
${\mathcal K}(X,Y) \subsetneqq {\mathcal L}(X,Y)$ implies that ${\mathcal K}(X,Y)$
is not complemented in ${\mathcal L}(X,Y)$. They prove that the answer is affirmative if
$Y$ has an unconditional basis.

A.E.~Tong \cite{T} proved that if $X$ has an unconditional basis, $Y$ is any Banach space,
and there exists a noncompact operator in ${\mathcal L}(X,Y^*)$ then $c_0$ embeds
in ${\mathcal L}(X,Y^*)$. This result was generalized by Kalton \cite{K} as follows:
Let $X$ be a Banach space with an unconditional finite dimensional expansion of the
identity (i.e. there is a sequence $(A_n) \subseteq {\mathcal L}(X)$ of finite rank operators
such that for every $x \in X$, $x= \sum_{n=1}^\infty A_n x$ unconditionally), and $Y$ be
any Banach space. If there is a noncompact operator in ${\mathcal L}(X,Y)$ then
$\ell_{\infty}$ embeds in ${\mathcal L}(X,Y)$. Moreover, the same assumptions imply that
${\mathcal K}(X,Y)$ is not complemented in ${\mathcal L}(X,Y)$.

G. Emmanuele \cite{E2} proved the following result: Assume that $X$ and $Y$ satisfy one of
the following two assumptions:
\begin{itemize}
\item[(i)] $X$ is a ${\mathcal L}_{\infty}$ space and $Y$ is a closed subspace of a
${\mathcal L}_1$ space, or
\item[(ii)] $X=C[0,1]$ and $Y$ is a space with cotype $2$.
\end{itemize}
If there is a noncompact operator in ${\mathcal L}(X,Y)$ then $\ell_\infty$
embeds in ${\mathcal L}(X,Y)$. Moreover, under the same assumptions
${\mathcal K}(X,Y)$ is not complemented in ${\mathcal L}(X,Y)$. Related results are also
in \cite{E1}.

%-------------------------------------New Section ------------------------------------

\section{Quasisubsymmetric Sequences}

In this section we introduce the notions, study, and give examples
of basic quasisubsymmetric sequences.  The main results are Theorems
$\ref{th:QSSEmbedding}$ and $\ref{thm:lastone}$ and Corollary
$\ref{th:lastone}$. Let $c_{00}$ denote the linear space of finitely
supported sequences.  If $(x_n)$ and $(y_n)$ are two basic sequences
in a Banach space and $C>0$, we say that $(x_n)$ $C$-dominates
$(y_n)$ if $\|\sum a_n y_n\| \leq C\|\sum a_nx_n\|$ for all $(a_n)
\in c_{00}$. We say that $(x_n)$ dominates $(y_n)$ if $(x_n)$
$C$-dominates  $(y_n)$ for some $C>0.$

\begin{defn} \label{defn:QSS}
A basic sequence $(e_n)_{n\in\mathbb{N}}$ in some Banach space $X$ is said to
be \textbf{quasisubsymmetric} if and only if
for any two increasing sequences $(k_n)$ and $(\ell_n)$  of
positive integers with $k_n \leq \ell_n$ for all $n$, we have that $(e_{k_n})$ dominates
$(e_{\ell_n})$, i.e. there exists a constant $S_{(k_n),(\ell_n)}>0$ such that
$$
\left\|\sum a_n e_{\ell_n}\right\|\leq S_{(k_n),(\ell_n)}\left\|\sum a_n e_{k_n}\right\|
\text{ for all }(a_n) \in c_{00}.
$$
\end{defn}

While in this definition the constant $S_{(k_n),(\ell_n)}$ depends
on the two subsequences, the dependence can be removed if $(e_n)$ is
seminormalized as the next result shows.

\begin{lem} Let $(e_n)$ be a seminormalized quasisubsymmetric basic sequence.  Then
there exists a constant $S>0$ such that for any two increasing
sequences $(k_n)$ and $(\ell_n)$ of positive integers with $k_n \leq
\ell_n$ for all $n$, we have that
\begin{equation}
\begin{split}
\left\|\sum a_n e_{\ell_n}\right\| \leq S\left\|\sum a_n
e_{k_n}\right\| \text{ for all }(a_n) \in c_{00}.
\label{eq:uniformS}
\end{split}
\end{equation}
\label{lem:uniformS}
\end{lem}

\begin{proof}
We will assume, via contradiction, that for all $S>0$ there exist
$(k_n)$ and $(\ell_n)$ increasing sequences of positive integers
with $k_n \leq \ell_n$ for all $n$, and there exists $(a_n)\in
c_{00}$ for which (\ref{eq:uniformS}) is not valid.

 \n {\bf Claim: } For every $S>0$ and every $N \in \N$ there exists $(a_n) \in c_{00}$ with $a_n =
0$ for all $n<N$ and there exist two increasing sequences $(k_n)$
and $(\ell_n)$ of positive integers with $k_n \leq \ell_n$ for all
$n$, such that
$$\left\|\sum a_n e_{\ell_n}\right\| > S\left\|\sum a_n e_{k_n}\right\|. $$

To prove this we fix $S$ and $N$ and let $\tilde S = S(1+B) +
(N-1)2BD/d,$ where $B$ is the basis constant of the basic sequence
$(e_n)$ and $d\leq \|e_n\| \leq D$ for all $n$. By assumption there
exists $(a_n) \in c_{00}$ and two increasing sequences $(k_n)$ and
$(\ell_n)$ of positive integers with $k_n \leq \ell_n$ for all $n$,
such that
$$\left\|\sum a_n e_{\ell_n}\right\| > \tilde S \left\|\sum a_n e_{k_n}\right\|.$$
We can assume without loss of generality that $\|\sum a_ne_{k_n}\| =
1.$  By the triangle inequality,
\begin{equation}
\begin{split}
\|\sum_{n=N}^\infty a_ne_{k_n} \|  &\leq 1+B.\label{eq:unc1}
\end{split}
\end{equation}

\n Also,
\begin{equation} \label{eq:coeff}
|a_n| \leq 2B/d  \mbox{ for all } n. 
\end{equation}

\n Furthermore,
\begin{equation*}
\begin{split}
\|\sum_{n=N}^\infty a_ne_{\ell_n} \|  &\geq \|\sum_{n=1}^\infty a_ne_{\ell_n} \|
             - \|\sum_{n=1}^{N-1} a_ne_{\ell_n} \|\\
    &> \tilde S \|\sum_{n=1}^\infty a_ne_{k_n} \| - (N-1)2BD/d \mbox{ (by (\ref{eq:coeff}))}\\
    &\geq S(1+B) + (N-1)2BD/d - (N-1)2BD/d \\
    &\geq S \|\sum_{n=N}^\infty a_ne_{k_n} \|\mbox{ (by ($\ref{eq:unc1}$))}
\end{split}
\end{equation*}

which finishes the proof of the claim.

 Let $N_1 = 1$ and $S_1 = 2$.  By the claim there exists an
$(a_n^{(1)}) \in c_{00}$ and two increasing sequences $(k_n^{(1)})$
and $(\ell_n^{(1)})$ of positive integers with $k_n^{(1)} \leq
\ell_n^{(1)}$ for all $n$, such that
$$\left\|\sum a_n^{(1)} e_{\ell_n^{(1)}}\right\| > 2\left\|\sum a_n^{(1)} e_{k_n^{(1)}}\right\| = 2.$$

Let $N_2 = \max (\ell_n^{(1)})$ where the maximum is taken over all
$n$ in the support of $(a_n^{(1)})$ and let $S_2 = 4$.  By the claim
there exists an $(a_n^{(2)}) \in c_{00}$ with $a_n^{(2)}=0$ for all
$n < N_2$  and two increasing sequences $(k_n^{(2)})$ and
$(\ell_n^{(2)})$ of positive integers with $k_n^{(2)} \leq
\ell_n^{(2)}$ for all $n$, such that
$$\left\|\sum a_n^{(2)} e_{\ell_n^{(2)}}\right\| > 4\left\|\sum a_n^{(2)} e_{k_n^{(2)}}\right\| = 4.$$

Continue in this manner with $N_m = \max (\ell_n^{(m-1)})$ where the
maximum is taken over all $n$ in the support of $(a_n^{(m-1)})$ and
let $S_m = 2^m$.

Let $a_n=a_n^{(m)}/2^m$, $k_n = k_n^{(m)}$ and $\ell_n =
\ell_n^{(m)}$ for all $n \in [N_m,N_{m+1})$ and for all $m \in
\mathbb{N}$.  Note that $\sum a_ne_{k_n}$ converges and has norm at
most equal to 2. However $\sum a_ne_{\ell_n}$ does not converge,
contradicting that $(e_n)$ is quasisubsymmetric.
\end{proof}

The authors do not know whether every seminormalized
quasisubsymmetric sequence can be equivalently renormed so that the
constant $S$ in the statement of Lemma 2.2 can be taken to be equal
to 1.  If this is true then the proof of Theorem
\ref{th:QSSEmbedding} can be slightly simplified.

Some examples of spaces with quasisubsymmetric basic sequences follow.

\begin{itemize}
\item[(a)] If any two subsequences of a basic sequence $(e_n)$ are naturally isomorphic
then $(e_n)$ is quasisubsymmetric; (the examples (b) and (e) below satisfy this condition).
In particular any subsymmetric sequence is quasisubsymmetric.
\item[(b)] There are conditional quasisubsymmetric sequences, like the summing basis.
\item[(c)] The sequence of the biorthogonal functionals of the basis of
Schreier's\footnote{An excellent resource for Schreier's space and Tsirelson's space is
\cite{CS}.} space is quasisubsymmetric.  We say that a finite subset $E$ of $\N$ is
a {\bf Schreier set} if the cardinality of $E$ is less than or equal to the minimum of $E$.
For $x=(x_k)  \in c_{00}$ we define the norm
$$\|x\| = \sup_{E~is~Schreier} \sum_{k \in E} |x_k|.$$
Schreier's space is the completion of $c_{00}$ equipped with this norm.

\item[(d)] The sequence of the biorthogonal functionals of the basis of
Tsirelson's\footnotemark[\value{footnote}] space \cite{FJ} is quasisubsymmetric.
There exists a norm $\| \cdot \|$ satisfying
$$\|x\| = \max\left \{\|x\|_{c_0}, \frac{1}{2} \sup\sum_{j=1}^k\|E_jx\|\right \}$$
for all $x \in c_{00}$, where $\| \cdot \|_{c_0}$ denotes the norm of $c_0$ and the supremum
 is taken over all sequences of sets $(E_j)_{j=1}^k$ such that $\max E_j < \min E_{j+1}$ and
$(\min E_j)_{j=1}^k$ is a Schreier set.  Also $E_jx$ denotes the natural projection of $x$ on
$E_j$.
Tsirelson's space is the completion of $c_{00}$ equipped with this norm.

\item[(e)] The James' space has a boundedly complete basis when the norm is defined as
$$\|x\| = \sup\left ( \sum_{k=0}^n \left (\sum_{i=p_k}^{p_{k+1}-1}
        x_i \right )^2 \right )^{1/2}$$
\n where $x=(x_i)$ and the supremum is taken over all positive integers n and all sequences
of integers
such that $0=p_0<p_1<\ldots<p_{n+1}$.  Since the basis is boundedly complete
\cite[page 50]{FG} $c_0 \not \hookrightarrow J$.  Also it is conditional
\cite[pages 49-50]{FG}.  Obviously every two subsequences of the basis are isometric.
Thus our Theorem~\ref{th:QSSEmbedding} gives an easy way to see
$\ell_\infty \hookrightarrow \mathcal{L}(J)$.
\end{itemize}

One of our main results of this section is Theorem~$\ref{th:QSSEmbedding}$ which gives
sufficient conditions on $X$ for $\ell_\infty \hookrightarrow \mathcal{L}(X).$  If, however,
we are only interested in $\mathcal{L}(X)$ being non separable, then the assumptions of
Theorem $\ref{th:QSSEmbedding}$ can be weakened as the next result shows.

\begin{prop}\label{prop:nseparable}
If a Banach space $X$ has a seminormalized basis which dominates its subsequences
then $\mathcal{L}(X)$ is non-separable.
\end{prop}

\begin{proof}
Let $(e_n)$ be a seminormalized basis for $X$ such that it dominates its subsequences and
let $C$ be its basis constant.  If $(k_n)$ is an increasing sequence in $\N$ then define
$T_{(k_n)}:X \rightarrow X$ by
$T_{(k_n)}(x) = \sum_{n=1}^\infty e_n^*(x) e_{k_n}$.  By our assumption each $T_{(k_n)}$ is
bounded, i.e. $T_{(k_n)} \in \mathcal L (X)$.  Let $(k_n)$, $(m_n)$
be two different increasing sequences in $\N$.  Then there exists some $j \in \N$ such that
$k_j \not = m_j$.  Thus
\begin{equation*}
\begin{split}
\|T_{(k_n)} - T_{(m_n)}\| & \geq  \left\|
\sum_{n=1}^\infty \left ( e_n^*\left(\frac{e_j}{\|e_j\|}\right) e_{k_n}-
e_n^*\left(\frac{e_j}{\|e_j\|}\right) e_{m_n} \right ) \right\| \\
 &= \frac{1}{\|e_j\|} \|e_{k_j} - e_{m_j}\| \geq
\frac{\|e_{k_j} - e_{m_j}\|}{\sup_n \|e_n\|} >  \frac{\inf_n\|e_n\|}{C\sup_n \|e_n\|} >  0.
\end{split}
\end{equation*}\end{proof}

We now present our first main result.

\begin{thm} \label{th:QSSEmbedding}
If a Banach space $X$ has a seminormalized quasisubsymmetric basis then
$\ell_{\infty} \hookrightarrow \mathcal{L}(X)$.
\end{thm}

In the proof of Theorem \ref{th:QSSEmbedding}, we will use the following remark 
due to Tong and Wilken \cite[Proposition 4]{TW} which is
mentioned in the introduction.

\begin{rem}  \label{rem:TongWilken}
Let $(w_n)$ be an unconditional basic sequence in a Banach space $X$, and assume that
there exists a noncompact operator in $\mathcal{L}(X,[(w_n)])$.  Then
$\ell_\infty \hookrightarrow \mathcal{L}(X)$.
\end{rem}

\begin{proof}[Proof of Theorem \ref{th:QSSEmbedding}] 
Let $(e_n)_{n \in \N}$ be a seminormalized quasisubsymmetric basis for $X$.
Since $(e_n)$ is bounded, by the $\ell_1$ theorem of H.P. Rosenthal \cite{R} there
exists a subsequence $(e_{k_n})_{n \in \N}$ which is either equivalent to the unit
vector basis of $\ell_1$ or it is weakly Cauchy.

In the first case, if $(e_{k_n})$ is K-equivalent to the unit vector basis of $\ell_1$, observe that
for all $(a_n) \in c_{00}$,

$$\frac{K}{S}\sum |a_n|\leq\frac{1}{S}
 \|\sum_{n=1}^\infty a_n e_{k_n} \|
\leq \| \sum_{n=1}^\infty  a_n e_n \| \leq \sup_n\|e_n\|\sum_{n=1}^\infty |a_n|$$

\n where $S$ is the constant provided by Lemma~$\ref{lem:uniformS}$.
Thus $(e_n)$ is equivalent to the unit vector basis of $\ell_1$, hence it is an
unconditional basis for $X$ and therefore $(e_n^* \otimes e_n)_{n=1}^\infty$ is equivalent
to the unit vector basis of $c_0$ (where for $x \in X$ and $x^* \in X^*$, $x^*\otimes x$
denotes the operator on $X$ defined by $(x^*\otimes x)(y) = x^*(y) x$ for every $y \in X$).

In the second case we assume that $(e_{k_n})_{n=1}^{\infty}$
is weakly Cauchy.  Define  the difference sequence $(z_n)$ of $(e_{k_n})$ by
$z_n = e_{k_{2n}}-e_{k_{2n-1}}$.
Then $(z_n)_{n\in \mathbb{N}}$ is weakly null.  Since $(e_n)$ is a seminormalized
basic sequence, we have that $(z_n)$ is also a seminormalized basic sequence,
say $0<d\leq\|z_n\|\leq D < \infty$ for all $n$. We claim that there exists
a subsequence of $(z_n)_{n\in \N}$ which is unconditional. In order to prove
this, we use an argument of A. Brunel and L. Sucheston \cite{BS}.

First by Mazur's theorem \cite[Corollary II.A.5]{W} there exist intervals of the
natural numbers $ I_1 < I_2 < ...$ (where $I_k<I_j$ means
$\max (I_k )< \min (I_j)$) and there exists $c_i \geq 0$ with
$\sum_{i\in I_j} c_i =1$ for all $j$, such that $\|\sum_{i\in I_j} c_i z_i \| \rightarrow 0$ as
$j \rightarrow \infty$.

By taking a subsequence and renaming we can assume without loss of generality, that

\begin{equation}
  \|\sum_{i\in I_n} c_i z_i \|<\frac{1}{2^n} . \label{eq:one}
\end{equation}

\noindent Let $(a_n) \in c_{00}$.  Let $F$ be a finite subsequence
of positive integers.  For $n \in \N$ let $m_n:= \min I_n$. We will show

\begin{equation} \label{eq:newstar}
\| \sum_{n \in \N \backslash F}a_n z_{m_n} \| \leq
(S + \frac{C}{d})\|\sum_{n=1}^\infty a_n z_{m_n}\|
        \mbox{ for all }a_n \in c_{00}
\end{equation}

\n where $C$ is the basis constant of $(z_n)$ and $S$ is the constant provided by
Lemma~$\ref{lem:uniformS}$ for the basic sequence $(e_n)$.

First note since $(e_n)$ is quasisubsymmetric that for all sequences
$(\ell_n)$ of positive integers with $\ell_n \in I_n$ (for all
$n$), and for all $(a_n) \in c_{00}$
$$ \|\sum_{n=1}^\infty a_n z_{\ell_n} \|\leq S \|\sum_{n=1}^\infty a_n z_{m_n} \|.$$

\n Let $(a_n) \in c_{00}$.  Assume that $F= \{t_1, t_2, \ldots, t_s\}$. For any
$(i_1, i_2, \ldots, i_s)$
with $i_n \in I_{t_n}$ for $n=1, \ldots , s$, define

$$V(i_{1}, i_{2}, \ldots, i_s) = \sum_{n \in \N \backslash F}a_nz_{m_n} +
        \sum_{n=1}^s a_{t_n}z_{i_n}$$
\n Thus if we sum over all $(i_1, i_2, \ldots, i_s)$ with $i_n \in I_{t_n}$ for
$n = 1, \ldots ,s$, we obtain

\begin{equation}
\sum c_{i_1}\cdots c_{i_s}V(i_1, \ldots, i_s) = \sum_{n \in \N \backslash F}a_nz_{m_n}
        + \sum_{n=1}^s a_{t_n}\sum_{i \in I_{t_n}}c_iz_i  \label{eq:star}
\end{equation}
\n Apply Lemma~\ref{lem:uniformS} for $k_n = m_n$ if $n \in \N$, $\ell_n=m_n$ if 
$n \in \N \backslash F$ and $\ell_n = i_n$ if $n \in F$, to obtain
\begin{equation*}
  \|V(i_{1}, \ldots, i_{s})\| \leq S\|\sum_{n=1}^\infty a_nz_{m_n}\| .
\end{equation*}

\n Thus

\begin{equation}
\begin{split}
\| \sum_{i_n \in I_{t_n} \text{ for }n = 1, \ldots , s} c_{i_1}\cdots c_{i_s} V(i_1, \ldots, i_s) \|
  &\leq
\sum_{i_n \in I_{t_n} \text{ for }n = 1, \ldots , s} c_{i_1}\cdots c_{i_s}\|V(i_1, \ldots, i_s) \|\\
  &\leq
\sum_{i_n \in I_{t_n} \text{ for }n = 1, \ldots , s} c_{i_1}\cdots c_{i_s} S
\|\sum_{n=1}^\infty a_n z_{m_n} \|\\
  &= S \|\sum_{n=1}^\infty a_n z_{m_n}\|. \label{eq:four}
\end{split}
\end{equation}

\n Also note

\begin{equation} \label{eq:three}
|a_n|d \leq \|a_nz_{m_n}\| = \|\sum_{i=1}^n a_iz_{m_i} - \sum_{i=1}^{n-1} a_iz_{m_i}\|
  \leq 2C\|\sum_{i=1}^\infty a_iz_{m_i}\| .
\end{equation}

\n Hence by ($\ref{eq:star}$),

\begin{equation}
\begin{split}
\| \sum_{i_n \in I_{t_n} \text{ for }n = 1, \ldots , s} c_{i_1} \cdots c_{i_s} V(i_1, \ldots, i_s) \|
&\geq \|\sum_{n \in \N \backslash F}a_nz_{m_n}\| - \sum_{n=1}^s|a_{t_n}| \,
            \| \sum_{i \in I_{t_n}} c_iz_i  \| \\
&\geq   \|\sum_{n \in \N \backslash F}a_nz_{m_n}\| -
  \frac{2C}{d} \|\sum_{i=1}^\infty a_i z_{m_i}\| \sum_{n=1}^s \frac{1}{2^{t_n}} \\
& \mbox{ (by (\ref{eq:one}) and (\ref{eq:three}))}\\
&\geq   \|\sum_{n \in \N \backslash F}a_nz_{m_n}\| -
  \frac{2C}{d}\|\sum_{n=1}^\infty  a_nz_{m_n}\|.
        \label{eq:five}
\end{split}
\end{equation}

\n So by ($\ref{eq:four}$) and ($\ref{eq:five}$) we obtain ($\ref{eq:newstar}$) which proves
that $(z_{m_n})$ is unconditional.

Let $P \in \mathcal{L}(X, [(z_{m_n})])$ be defined by
$$P\left(\sum_{n=1}^{\infty} a_n e_n\right)=\sum_{n=1}^{\infty} a_n z_{m_n}.$$
\n Notice that $P$ is bounded since
\begin{align*}  \left\|P \left(\sum_{n=1}^{\infty} a_n e_n \right) \right\|
        &\leq
\|a_1 e_{k_{2m_1}} + a_2 e_{k_{2m_2}} + \cdots \| +
\|a_1 e_{k_{2m_1-1}} + a_2 e_{k_{m_2 -1}} + \cdots \| \\
        &\leq 2S\| \sum_{n=1}^{\infty} a_n e_n \|.
\end{align*}

\n Since $P$ is a noncompact operator and $(z_{m_n})$ is unconditional,
by Remark $\ref{rem:TongWilken}$, we have
$\ell_\infty \hookrightarrow \mathcal{L}(X).$\end{proof}

By observing that for a Banach space $X$, the map from ${\mathcal L}(X)$ to 
${\mathcal L}(X^*)$ given by $T\mapsto T^*$, where $T^*$ is the adjoint operator, 
is an isometric embedding, 
Theorem \ref{th:QSSEmbedding} gives the following corollary.

\begin{cor} \label{cor:LastCor}
If a  Banach space  $X$ has a seminormalized quasisubsymmetric basis
then $\ell_\infty \hookrightarrow \mathcal{L}(X^*)$.
\end{cor}

The proof of Theorem $\ref{th:QSSEmbedding}$ gives the following two corollaries:

\begin{cor}
If $(e_n)_{n\in \N}$ is a seminormalized quasisubsymmetric basic sequence then either
$(e_n)_{n\in \N}$ is equivalent to the unit vector basis of $\ell_1$ or there exists a
subsequence $(e_{k_n})$ such that $(e_{k_{2n}}-e_{k_{2n-1}})_{n\in \N}$ is unconditional.
\label{cor:UncondSubseq}
\end{cor}

\begin{cor}
If $(e_n)_{n\in \N}$ is a seminormalized quasisubsymmetric weakly
null basic sequence then there exists a subsequence $(e_{k_n})$ of
$(e_{n})$ which is unconditional.
\end{cor}

\begin{thm}
If a Banach space $X$ has a seminormalized basis $(e_n)$ and the sequence $(e_n^*)$
of the biorthogonal functionals is quasisubsymmetric then
$\ell_{\infty} \hookrightarrow \mathcal{L}(X)$. \label{thm:lastone}
\end{thm}

\begin{proof}
By results of Kalton \cite{K} mentioned in the introduction, it suffices to show that
$c_0 \hookrightarrow \mathcal{L}(X)$.  By Corollary \ref{cor:UncondSubseq} we separate
two cases:

If $(e^*_n)_{n\in\N}$ is equivalent to the unit vector basis of $\ell_1$, we
know that $X=[(e_n)_{n\in\N}]$ is isomorphic to $c_0$.
Therefore $c_0\hookrightarrow \mathcal{L}(X)$.

If $(e^*_{k_{2n}}-e^*_{k_{2n-1}})_{n \in \N}$ is unconditional, let
$K$ be its unconditionality constant and proceed by showing that
$((e^*_{k_{2n}}-e^*_{k_{2n-1}}) \otimes e_n)_{n\in \N}$ is
equivalent to the unit vector basis of $c_0$:  For $(a_n)_{n\in \N}
\in c_{00}$,

\begin{equation*}
\begin{split}
  \left \| \sum a_n (e^*_{k_{2n}}-e^*_{k_{2n-1}})\otimes e_n \right\|
 & = \left \|\left( \sum a_n (e^*_{k_{2n}}-e^*_{k_{2n-1}})\otimes e_n \right)^*\right\| \\
 & = \left \| \sum a_n e_n\otimes (e^*_{k_{2n}}-e^*_{k_{2k-1}}) \right\|
        \mbox{(where $e_n$ is considered in $X^{**}$)}\\
 & \leq 2\left\| \sum a_n e_n (x^*) (e^*_{k_{2n}}-e^*_{k_{2n-1}}) \right\|
        \mbox{(for some $x^* \in [(e_n^*)]$, $\|x^*\| \leq 1$)}\\
 & \leq 2K\max|a_n| \left \| \sum e_n (x^*) (e^*_{k_{2n}}-e^*_{k_{2n-1}}) \right \| \\
 & \leq 2K\max|a_n| \left(\left \| \sum e_n (x^*) e^*_{k_{2n}} \right \| +
\left \| \sum e_n (x^*) e^*_{k_{2n-1}} \right \| \right) \\
 & \leq 2K\max|a_n| \left(S_{(n),(k_{2n})} \left \| \sum e_n (x^*) e^*_n \right \| +
S_{(n),(k_{2n-1})} \left \| \sum e_n (x^*) e^*_n \right \| \right) \\
 & =2 K (S_{(n),(k_{2n})} +S_{(n),(k_{2n-1})} ) \| x^* \|  \max|a_n|\\
&\leq 2 K (S_{(n),(k_{2n})} +S_{(n),(k_{2n-1})} )\max|a_n|.
\end{split}
\end{equation*}

Verifying the reverse inequality (with a different constant) is trivial and both inequalities
together give us that $c_0 \hookrightarrow \mathcal{L}(X)$.
\end{proof}

\begin{cor}
If a Banach space $X$ has a seminormalized quasisubsymmetric basis $(e_n)$ then
$\ell_{\infty} \hookrightarrow \mathcal{L}([(e_n^*)])$, where $(e_n^*)$ is the sequence
of the biorthogonal functionals to $(e_n)$. \label{th:lastone}
\end{cor}

\begin{proof}Consider $(e_n)$ in $X^{**}$: it is seminormalized, quasisubsymmetric and 
biorthogonal to $(e_n^*)$.  Therefore by Theorem \ref{thm:lastone} we have 
$\ell_{\infty} \hookrightarrow \mathcal{L}([(e_n^*)])$.
\end{proof}

Finally we make some remarks related to the question of whether $\ell_1$ 
embeds into the space of nuclear operators of a Banach space $X$.
If $X$ is a Banach space let $(\mathcal{N}(X),\nu(\cdot))$ denote
the space of nuclear operators on $X$, namely

$$\mathcal{N}(X)=\{T\in\mathcal{L}(X):T=\sum_{n=1}^{\infty}x^*_n\otimes y_n ~\mbox{where}~ x^*_n\in X^*, ~
y_n\in X ~\mbox{and}~ \sum_{n=1}^{\infty}\|x^*_n\|\|y_n\| < \infty
\}$$
$$\nu(T)=\inf\{\sum_{n=1}^{\infty}\|x^*_n\|\|y_n\|: \sum_{n=1}^{\infty}x^*_n\otimes y_n \mbox{ is
a representation of }T\}.$$

The following remark is well known:

\begin{rem}
If $X$ is a Banach space with the Approximation Property then
$\mathcal{N}(X)^*$ is isometric to $\mathcal{L}(X^*)$. Moreover, if
$T \in \mathcal{L}(X^*)$ then the action of T as a functional on
$\mathcal{N}(X)$ gives the trace of $TS^*$, $\text{tr}\, (TS^*)$,
when it is applied to $S \in \mathcal{N}(X)$. \label{th:nandl}
\end{rem}

Indeed, it is shown in \cite[Proposition 4.6 (i) $\Leftrightarrow$
(ii), cf. Corollary 4.8 (a)]{Ry} that if the Banach space $X$ has
the Approximation Property then ${\mathcal N}(X)$ is naturaly
isometric to the projective tensor product $X^* \hat \otimes X$.
Also in \cite[page 24]{Ry} it is shown that for any Banach space
$Y$, $( Y^* \hat \otimes Y)^*$ is isometric to ${\mathcal L}(Y^*)$,
where for $T \in {\mathcal L}(Y^*)$, the action of $T$ as a
functional on $\sum_{i=1}^n y_i^* \otimes y_i \in Y^* \hat \otimes
Y$ gives $\sum_{i=1}^n (Ty_i^*)(y_i)$.

By Remark \ref{th:nandl} we have that for a Banach space $X$ with the 
Approximation Property, $\ell_1$ embeds complementably in 
${\mathcal N}(X)$ if and only if $\ell_\infty$ embeds in ${\mathcal L}(X^*)$.

The next remarks give sufficient conditions on a Banach space $X$ which 
ensure that $\ell_1$ embeds complementably into $\mathcal{N}(X)$.  A similar result is
provided by Tong \cite[Theorem 1.5]{T} who proved that for a dual
Banach space $X^*$ with an unconditional basis, there exists a
noncompact operator in $\mathcal{L}(X^*)$ if and only if $\ell_1$
embeds in $X^* \hat \otimes X$ complementably.  

\begin{rem}
Let $X$ be a Banach space which has the Approximation Property and assume that
$\ell_\infty$ embeds in ${\mathcal L}(X)$.  Then
$\ell_1 \hookrightarrow \mathcal{N}(X)$ complementably.
\label{th:last}
\end{rem}

\begin{proof}
Since the map from ${\mathcal L}(X)$ to ${\mathcal L}(X^*)$ given by 
$T \mapsto T^*$ is an isometric embedding, we have that  
$\ell_\infty$ embeds in ${\mathcal L}(X^*)$. Hence by \cite{BP} and Remark~\ref{th:nandl}
we obtain that $\ell_1$ embeds complementably in ${\mathcal N}(X)$.
\end{proof}

\begin{cor}  \label{th:MainNuc}
If a Banach space $X$ has a seminormalized quasisubsymmetric
basis then $\ell_1 \hookrightarrow \mathcal{N}(X)$ complementably.
\end{cor}

\begin{proof} 
Combine Theorem~\ref{th:QSSEmbedding} and Remark~\ref{th:last}.
\end{proof}

\begin{rem} \label{new}
Let $X$ be a Banach space which has the Approximation Property and contains
an unconditional basic sequence $(e_n)$. Suppose that there exists a noncompact
operator in ${\mathcal L}(X,[(e_n)])$. 
 Then
$\ell_1 \hookrightarrow \mathcal{N}(X)$ complementably.
\end{rem}

\begin{proof}
By Remark \ref{rem:TongWilken} we have that $\ell_\infty$ embeds in 
${\mathcal L}(X)$. Now Remark~\ref{th:last} finishes the proof.
\end{proof}

{\footnotesize
\noindent
Department of Mathematics, University of South Carolina, Columbia, SC 29208, \\
\noindent giorgis@math.sc.edu, beanland@math.sc.edu, dilworth@math.sc.edu, sanacory@math.sc.edu
}

\end{document}